\documentclass[10pt]{article}
\usepackage{amsfonts}
\usepackage{latexsym}
\usepackage{cite}
\usepackage{amsmath,amsfonts,latexsym,amssymb}
\usepackage[mathscr]{eucal}
\usepackage{cases,color}
\usepackage{amsthm}

\usepackage[bf,small]{caption2}
\usepackage{float}
\usepackage{graphicx}
\usepackage{amsmath}
\usepackage{amssymb}
\usepackage[all]{xy}

\newtheorem{theorem}{theorem}[section]

\newtheorem{lem}[theorem]{Lemma}

\newtheorem{rmk}[theorem]{Remark}
\newtheorem{thm}[theorem]{Theorem}

\begin{document}

\title{\vspace{-2cm}\textbf{On the structure of Kauffman bracket skein algebra of a surface}}
\author{\Large Haimiao Chen}
\date{}
\maketitle

\begin{abstract}
  Suppose $R$ is a commutative ring with identity and a fixed invertible element $q^{\frac{1}{2}}$ such that $q+q^{-1}$ is invertible.
  For an oriented surface $\Sigma$, let $\mathcal{S}(\Sigma;R)$ denote the Kauffman bracket skein algebra of $\Sigma$ over $R$. It is shown that to each embedded graph $G\subset\Sigma$ satisfying that $\Sigma\setminus G$ is homeomorphic to a disk and some other mild conditions, one can associate a generating set for $\mathcal{S}(\Sigma;R)$, and the ideal of defining relations is generated by relations of degree at most $6$ supported by certain small subsurfaces.

\end{abstract}

\section{Introduction}


Let $R$ be a commutative ring with identity and a fixed invertible element $q^{\frac{1}{2}}$. For an oriented surface $\Sigma$, the {\it Kauffman bracket skein algebra} of $\Sigma$ over $R$, denoted by $\mathcal{S}(\Sigma;R)$, is defined as the quotient of the free $R$-module generated by isotopy classes of framed links embedded in $\Sigma\times(0,1)$ by the submodule generated by the skein relations
\begin{figure}[h]
  \centering
  \includegraphics[width=8.5cm]{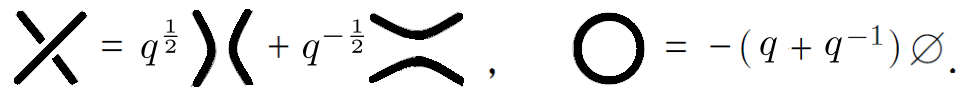}
\end{figure}
\\
As a convention, $R$ is identified with $R{\O}\subset\mathcal{S}(M;R)$ via $\eta\mapsto\eta{\O}$. Each link is equipped with the blackboard framing.

Given $L_1,L_2\subset\Sigma\times(0,1)$, their product $L_1L_2$ is defined by stacking $L_1$ over $L_2$ in the $[0,1]$ direction. Using the local relations, each element of $\mathcal{S}(\Sigma;R)$ can be written as a $R$-linear combination of multi-curves in $\Sigma$, where a multicurve is considered as a link in $\Sigma\times(0,1)$ by identifying $\Sigma$ with $\Sigma\times\{1/2\}$.

The description of the structure of $\mathcal{S}(\Sigma;\mathbb{Z}[q^{\pm\frac{1}{2}}])$ was raised as \cite{Ki97} Problem 1.92 (J) (proposed by Bullock and Przytycki), and also \cite{Oh02} Problem 4.5. A finite set of generators was given by Bullock \cite{Bu99}. So the problem is really to determine the defining relations.

Let $\Sigma_{g,h}$ denote the oriented surface of genus $g$ with $h$ boundary circles; let $\Sigma_g=\Sigma_{g,0}$.
In the previous work \cite{Ch22}, we found a generating set of $\mathcal{S}(\Sigma_{0,n+1};R)$ consisting of $n+\binom{n}{2}+\binom{n}{3}$ elements,
and showed that relations among them can be ``localized", in the sense that the ideal of defining relations is generated by certain relations of degree at most $6$. In this paper, we extend the result to $\Sigma_{g,h}$ for all $g,h$. 
Setting $q^{\frac{1}{2}}=-1$ will lead to a presentation for $\mathbb{C}[\mathcal{X}_{{\rm SL}(2,\mathbb{C})}(\Sigma_{g})]$,
which was regarded as a difficult problem and no result is seen in the literature.

When $q+q^{-1}$ is not assumed to be invertible, in principle a presentation for $\mathcal{S}(\Sigma_{g,h};\mathbb{Z}[q^{\pm\frac{1}{2}}])$ can still be obtained without essential difficulty.
Through this approach, the structure of $\mathcal{S}(\Sigma_{0,5};\mathbb{Z}[q^{\pm\frac{1}{2}}])$ had been determined in \cite{Ch23}.

The content is organized as follows. In Section 2 we introduce some necessary notations and conventions (most of which have been given in \cite{Ch22}), and propose the notion ``cutting system"; we define admissible expressions for all connected $1$-submanifolds in $\Sigma\times[0,1]$ for all surface $\Sigma$. In Section 3 we formulate Theorem \ref{thm:relation} which (roughly) asserts that relations can be ``localized", and outline the process of proving it; the core of the proof is a ``tortuous induction" which was implemented in \cite{Ch22}. In Section 4, several examples of cutting systems on surfaces are given. 
Section 5 is an appendix, containing detailed formulas supporting Lemma \ref{lem:substitution}.

\section{Generating set arising from a cutting system}

\subsection{Preparation}

Throughout, denote $q^{-1}$ as $\overline{q}$ (and also denote $q^{-\frac{1}{2}}$ as $\overline{q}^{\frac{1}{2}}$, etc), to make expressions as compact as possible. 
Let $\Sigma=\Sigma_{g,h}$. Let $\pi:\Sigma\times[0,1]\to\Sigma$ denote the projection.
We usually abbreviate $\mathcal{S}(\Sigma_{g,h};R)$ to $\mathcal{S}$.

For a finite set $Y$, let $\#Y$ denote its cardinality.

We always use `isotopy' to mean ordinary isotopy, without further explanations.
A {\it fine isotopy} is a family of homeomorphisms $\varphi_t:\Sigma\times[0,1]\to\Sigma\times[0,1]$, $0\le t\le 1$, such that $\varphi_0={\rm id}$, $\pi(\varphi_t(\mathsf{a},z))=\pi(\varphi_t(\mathsf{a},0))$ for all $(\mathsf{a},z)\in\Sigma\times[0,1]$, and $\varphi_t(\Gamma_k)=\Gamma_k$ for all $t,k$. This notion is introduced for keeping track of the information on crossings.
Given $1$-submanifolds $X,X'\subset\Sigma\times[0,1]$, call them {\it finely isotopic} and denote $X\approx X'$, if there exists a fine isotopy $\varphi_t$ with $\varphi_1(X)=X'$; in this case, $\lambda(X)=\lambda(X')$.

Call a (most likely non-connected) embedded graph $G\subset\Sigma$ a {\it cutting system} for $\Sigma$ if the following conditions hold: (i) $\Sigma\setminus G$ is homeomorphic to a disk; 
(ii) each point in $G\cap\partial\Sigma$ is a univalent vertex; (iii) the arcs incident to each interior vertex of $G$ (i.e. one lying in the interior of $\Sigma$) belong to $2$, $3$, or $4$ edges; (iv) each edge of $G$ is oriented.

Fix a cutting system $G$.
Numerate the edges of $G$ as $\gamma_1,\ldots,\gamma_n$.
Let $\Gamma_i=\gamma_i\times[0,1]\subset\Sigma\times[0,1]$; let $\Gamma=\bigcup_{i=1}^n\Gamma_i$.

Given $\{i_1,\ldots,i_r\}\subset\{1,\ldots,n\}$, let $\Sigma(i_1,\ldots,i_r)$ denote the surface obtained from cutting $\Sigma$ along the $\gamma_j$'s for $j\in\{1,\ldots,n\}\setminus\{i_1,\ldots,i_r\}$.

A $1$-submanifold $X\subset\Sigma\times[0,1]$ is always assumed to be {\it in generic position}, meaning that $\pi(X)$ is stable under small perturbations of $X$. In particular, $\pi(X)\cap G=\emptyset$.
Let ${\rm Cr}(X)$ denote the set of crossings of $X$. Let $|X|=\#(X\cap\Gamma)$, called the {\it degree} of $X$,
and let ${\rm cn}(X)=\#{\rm Cr}(X)$. Call $X$ {\it simple} if ${\rm cn}(X)=0$. Define $\lambda(X)=(|X|,{\rm cn}(X))$.

Each arc is assumed to be oriented, although sometimes the orientation is irrelevant. For an arc $F$, write $\partial F=\{\partial_-F,\partial_+F\}$ so that $F$ is oriented from $\partial_-F$ to $\partial_+F$. Let $\overline{F}$ denote the arc obtained from $F$ by reversing its orientation.

Given $\Omega=\sum_ia_iX_i$ with $0\ne a_i\in R$ and $X_i$ being a $1$-submanifold, let ${\rm md}_{\Omega}(v)=\max_i{\rm md}_{X_i}(v)$, and let $|\Omega|=\sum_{v=1}^n{\rm md}_{\Omega}(v)$.
For linear combinations $\Omega,\Omega'$, write ${\rm md}_{\Omega}\le {\rm md}_{\Omega'}$ if ${\rm md}_{\Omega}(v)\le{\rm md}_{\Omega'}(v)$ for all $v$.

Introduce a linear order $\preceq$ on $\mathbb{N}^2$, by declaring $(m',c')\preceq(m,c)$ if either $m'<m$, or $m'=m$, $c'\le c$. Denote $(m',c')\prec(m,c)$ if $(m',c')\preceq(m,c)$ and $(m',c')\ne(m,c)$. If $\lambda(X)\prec\lambda(X')$, then we say that $X$ is {\it simpler} than $X'$.

Let $\mathcal{F}$ denote the free $R$-algebra generated by fine isotopy classes of links, with multiplication defined via stacking.
Let $\theta':\mathcal{F}\twoheadrightarrow\mathcal{S}$ denote the canonical quotient. For $\mathfrak{u}_1,\mathfrak{u}_2\in\mathcal{F}$, we say $\mathfrak{u}_1=\mathfrak{u}_2$ in $\mathcal{S}$ if $\theta'(\mathfrak{u}_1)=\theta'_n(\mathfrak{u}_2)$.

Let $\mathcal{V}$ denote the free $R$-module generated by isotopy classes of multi-curves.
For a link $L$, let $\Theta(L)\in\mathcal{V}$ denote the linear combination of isotopy classes of multi-curves resulting from resolving all crossings of $L$. Extend this via linearity to $\Theta:\mathcal{S}\to\mathcal{V}$; by definition, it descends to $\Theta:\mathcal{S}\to\mathcal{V}$. By \cite{CM12,SW07}, $\Theta:\mathcal{S}\cong\mathcal{V}$ is an isomorphism of $R$-modules.



When $A,B$ are arcs with $\partial_+A=\partial_-B$ and $\partial_-A\ne \partial_+B$, let $AB$ denote the arc obtained by identifying $\partial_+A$ with $\partial_-B$, i.e. $AB=A\cup B$; when $\partial_+A=\partial_-B$ and $\partial_-A=\partial_+B$,
let ${\rm tr}(A B)=A\cup B$, which is a knot. In the construction of such kind, we may perturb $A$ or $B$ if necessary, to ensure $A\cup B$ to be generic. This convention will be always adopted.





Let $\mathbf{w}(F)$ denote the word in $\mathbf{x}_1^{\pm1},\ldots,\mathbf{x}_n^{\pm1}$ determined as follows: starting from $\partial_-F$, walk along $F$ guided by the orientation of $F$, record $\mathbf{x}_i$ (resp. $\mathbf{x}_i^{-1}$) whenever passing through $\Gamma_i$ from left to right (resp. from right to left), and multiply all the recorded $\mathbf{x}_i^{\pm1}$ together when reaching $\partial_+F$. Call $F$ {\it reduced} if $\mathbf{w}(F)$ is a reduced word; otherwise, call $F$ {\it reducible}.

Given simple curves $S,S'$, if there exist simple arcs $A,F,F'$ such that $|F|,|F'|\le 4$, $S={\rm tr}(AF)$, $S'={\rm tr}(AF')$, and ${\rm tr}(\overline{F}F')$ is a short simple curve enclosing some interior vertex of $G$, then we say that $S'$ is obtained from $S$ by sliding $F'$ to $F$ across the vertex.  

\begin{lem}\label{lem:isotopy}
Two simple curves are isotopic if and only if they can be related by fine isotopy, shrinking degree $2$ arcs, and sliding across interior vertices.
\end{lem}

\begin{proof}
The ``if" part is trivial.

For the ``only if" part, suppose $S,S'$ are isotopic simple curves. Remove a small disk around each interior vertex. Then $\Sigma$ becomes a surface with boundary $\Sigma^\circ$, and $G$ becomes a cutting system $G^\circ$ for $\Sigma^\circ$ which is a disjoint union of arcs. We may assume $S,S'\subset\Sigma^\circ$.

By sliding across interior vertices and shrinking degree $2$ arcs whenever necessary, we may assume that $S,S'$ are reduced and isotopic in $\Sigma^\circ$.

Then the proof of \cite{Ch22} Lemma 2.2 can be easily extended to $\Sigma^\circ$ here, to show that $S,S'$ are finely isotopic. 
\end{proof}

\subsection{Admissible expressions}


Let $\Delta$ denote the closure of $\Sigma\setminus G$, which by hyperthesis is a closed disk. Let ${\rm gl}:\Delta\to\Sigma$ stand for the gluing map.

Each simple arc $F\subset\Sigma$ of degree $3$ can be presented in the following way. Suppose $F$ intersects $\gamma_{i_1}$, $\gamma_{i_2}$, $\gamma_{i_3}$ consecutively in the order given by the orientation, with $i_1,i_2,i_3$ not necessarily distinct.
Take a homeomorphism $f:\Delta\to D^2$ (the unit disk) such that the $6$ points constituting $f({\rm gl}^{-1}(F\cap G))$ is equidistributed on $\partial D^2$.
Let $\widetilde{F}={\rm gl}^{-1}(F)$. Call $f(\widetilde{F})$ the {\it symbol} of $F$.

\begin{figure}
  \centering
  \includegraphics[width=12cm]{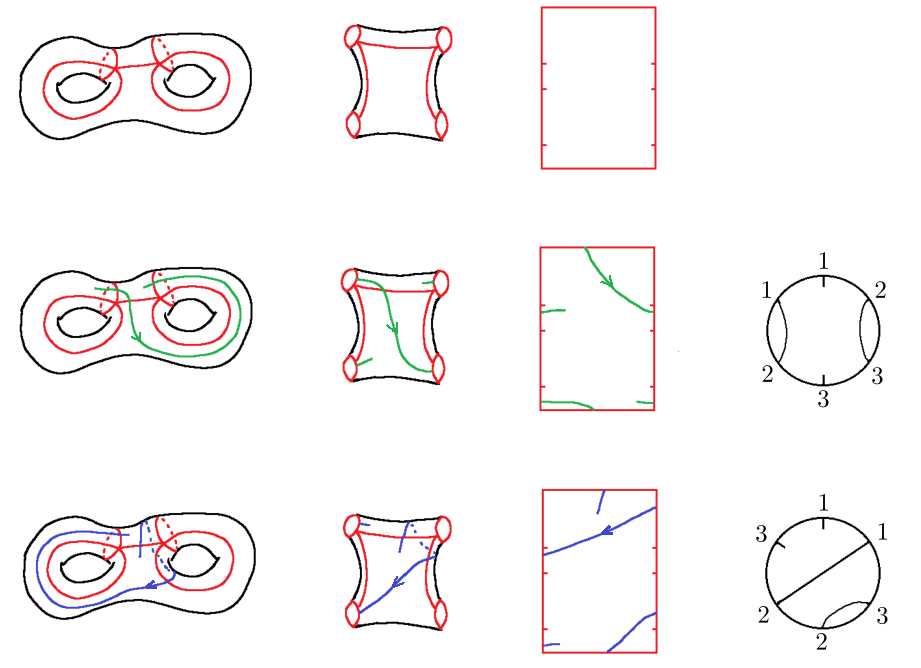}\\
  \caption{Upper row (from left to right): a cutting system of $\Sigma_2$; the $4$-holed sphere obtained from cutting two curves; the disk obtained from cutting along the curves.
  Middle and lower rows: degree $3$ simple arcs, with orientations indicated by arrows.}\label{fig:example}
\end{figure}

Two examples of arcs are shown in Figure \ref{fig:example}, where $\Sigma=\Sigma_2$.

\begin{figure}[H]
  \centering
  \includegraphics[width=12.5cm]{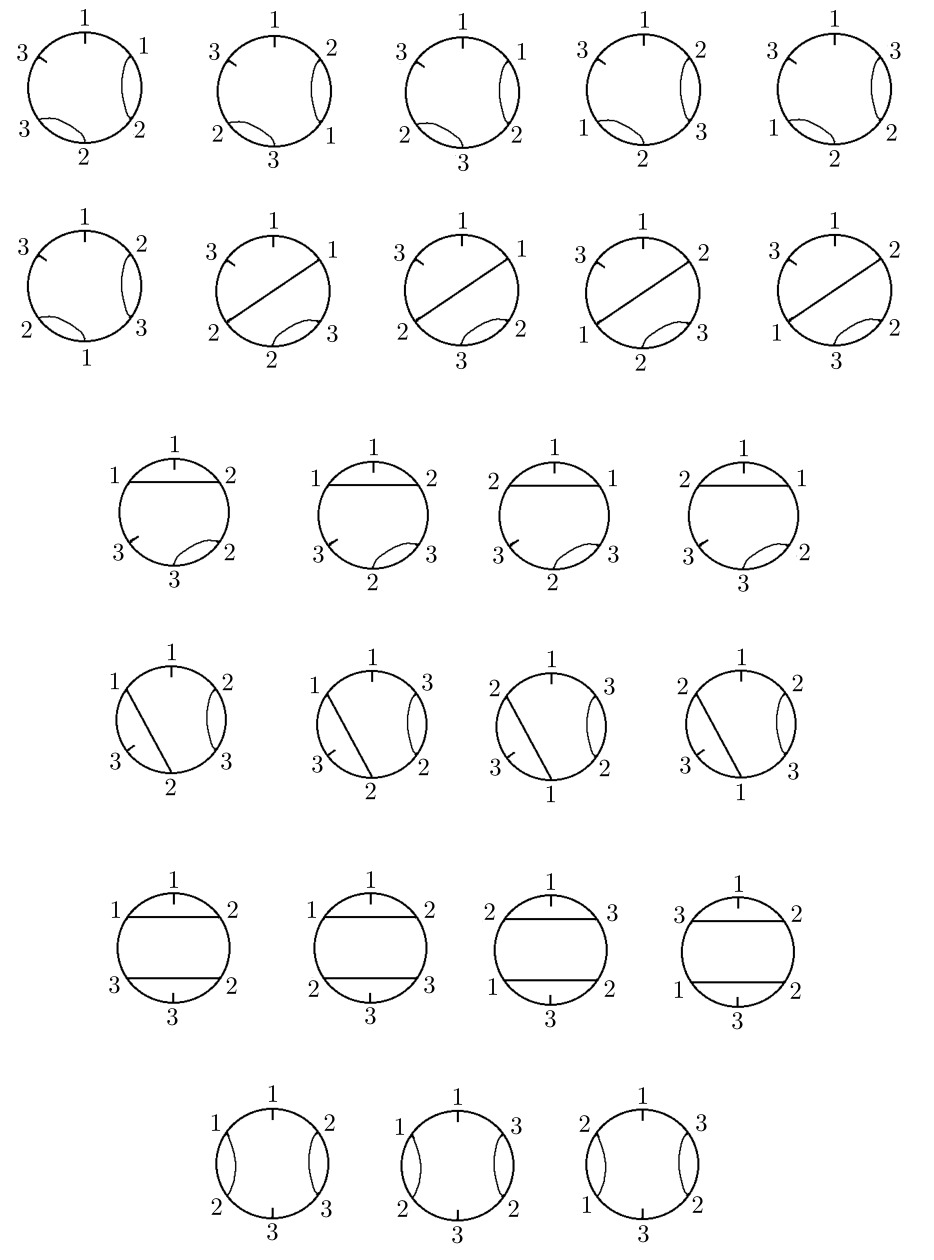}
  \caption{Various symbols for degree $3$ simple arcs. Denote the ones in the first two rows by ${\rm I}_1$ to ${\rm I}_{10}$; denote the ones in the third and fourth rows by ${\rm II}_1$ to ${\rm II}_8$; denote the ones in the last two rows by ${\rm III}_1$ to ${\rm III}_7$.
  }\label{fig:cases}
\end{figure}

Up to orientation-reversion of $F$ and self-homeomorphism of $D^2$, all possible symbols are shown in Figure \ref{fig:cases}.

Choose a point $\mathsf{x}_i\in\gamma_i$ for each $i$. Suppose ${\rm gl}^{-1}(\mathsf{x}_i)=\{\tilde{\mathsf{x}}_i^+,\tilde{\mathsf{x}}_i^-\}$.

For each $i$, take an arc $\tilde{t}_i\subset\Delta$ connecting $\tilde{\mathsf{x}}_i^+$ to $\tilde{\mathsf{x}}_i^-$, and let $t_i={\rm gl}(\tilde{t}_i)$.
For $i<j$, choose once for all a disjoint union $\tilde{t}_{ij}$ of two arcs connecting $\tilde{\mathsf{x}}_i^{\pm}$ to $\tilde{\mathsf{x}}_j^{\pm\epsilon}$, with $\epsilon\in\{\pm\}$, and let $t_{ij}={\rm gl}(\tilde{t}_{ij})$. For $i<j<k$, choose once for all
a disjoint union $\tilde{t}_{ijk}$ of three arcs connecting $\tilde{\mathsf{x}}_i^+$ to $\tilde{\mathsf{x}}_j^{\epsilon}$, connecting $\tilde{\mathsf{x}}_j^{-\epsilon}$ to $\tilde{\mathsf{x}}_k^{\o}$, and connecting $\tilde{\mathsf{x}}_k^{-\o}$ to $\tilde{\mathsf{x}}_i^-$, where $\epsilon,\o\in\{\pm\}$, and let $t_{ijk}={\rm gl}(\tilde{t}_{ijk})$.

\begin{rmk}
\rm When $\tilde{\mathsf{x}}_i^{\pm},\tilde{\mathsf{x}}_j^{\pm}$ are not interlaced, $\tilde{t}_{ij}$ is unique. Otherwise, there are two choices for $\tilde{t}_{ij}$; however, denoting the other choice by $\tilde{t}'_{ij}$, we have that ${\rm gl}(\tilde{t}'_{ij})$ equals a linear combination of $t_it_j$ and $t_{ij}$.

The situation for $\tilde{\mathsf{x}}_i^{\pm},\tilde{\mathsf{x}}_j^{\pm},\tilde{\mathsf{x}}_k^{\pm}$ is similar, although more complicated. It is not difficult to see that the ${\rm gl}(\tilde{t}'_{ijk})$ for any other choice $\tilde{t}'_{ijk}$ for $\tilde{t}_{ijk}$ can be written as a linear combination of $t_it_jt_k$, $t_it_{jk}$, $t_jt_{ik}$, $t_kt_{ij}$ and $t_{ijk}$.
\end{rmk}

Let $\mathcal{T}$ be the free $R$-algebra generated by $\mathfrak{T}$, where
$$\mathfrak{T}=\{t_{i_1\cdots i_r}\colon 1\le i_1<\cdots<i_r\le n,\ 1\le r\le 3\}.$$
Let $\theta:\mathcal{T}\to\mathcal{S}$ denote the canonical map.

Given $\mathsf{x},\mathsf{y}\in\Sigma\times\{z\}$ with $z\in\{0,1\}$, let $\mathcal{H}_z(\mathsf{x},\mathsf{y})$ denote the set of $1$-submanifolds
$X=C\cup L\subset\Sigma\times[0,1]$, where $L$ is a link and $C$ is an arc with $\partial_-C=\mathsf{x}$, $\partial_+C=\mathsf{y}$.
Let $\mathcal{S}(\mathsf{x},\mathsf{y})$ denote the $R$-module generated by relative isotopy classes of elements of $\mathcal{H}_z(\mathsf{x},\mathsf{y})$, modulo skein relations. Here a relative isotopy means an isotopy $\varphi_t$ of $\Sigma\times[0,1]$ with $\varphi_t(\mathsf{x})=\mathsf{x}$ and $\varphi_t(\mathsf{y})=\mathsf{y}$ for all $t$. Let $[X]\in\mathcal{S}(\mathsf{x},\mathsf{y})$ denote the element represented by $X$.
When $z=0$ (resp. $z=1$), $\mathcal{S}(\mathsf{x},\mathsf{y})$ is a left (resp. right) $\mathcal{S}$-module.
Let $\mathcal{M}(\mathsf{x},\mathsf{y})$ denote the subset of $\mathcal{H}_z(\mathsf{x},\mathsf{y})$ consisting of simple reduced arcs $C$ with $|C|\le 2$.

\begin{lem}\label{lem:substitution}
Let $F\in\mathcal{H}_0(\mathsf{x},\mathsf{y})$ be a degree $3$ arc.
\begin{enumerate}
  \item[\rm(i)] If $z=0$, then there exist $\mathfrak{a}_s\in\mathcal{T}$, $C_s\in\mathcal{M}(\mathsf{x},\mathsf{y})$ such that $[F]=\mathfrak{s}_u(F):=\sum_s\mathfrak{a}_s[C_s]$ in $\mathcal{S}(\mathsf{x},\mathsf{y})$, and ${\rm md}_{\mathfrak{s}_u(F)}\le {\rm md}_F$.
  \item[\rm(ii)] If $z=1$, then there exist $\mathfrak{b}_t\in\mathcal{T}$, $D_t\in\mathcal{M}(\mathsf{x},\mathsf{y})$ such that $[F]=\mathfrak{s}_d(F):=\sum_t[D_t]\mathfrak{b}_t$ in $\mathcal{S}(\mathsf{x},\mathsf{y})$, and ${\rm md}_{\mathfrak{s}_d(F)}\le {\rm md}_F$.
\end{enumerate}
\end{lem}

\begin{proof}
The result is an extension of \cite{Ch22} Lemma 3.1.

Say that ${\rm I}_1$ can be chopped if an arc of type ${\rm I}_1$ satisfies the assertion, and so forth.
As a key ingredient of the proof of \cite{Ch22} Lemma 3.1, it was shown that ${\rm I}_1, {\rm II}_1, {\rm III}_1$ can be chopped.
\begin{figure}[H]
  \centering
  \includegraphics[width=9cm]{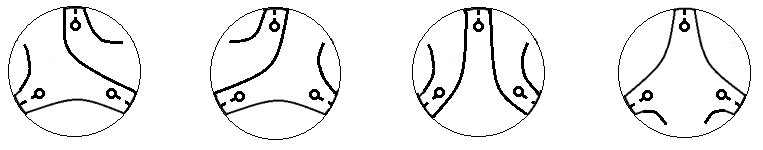}\\
  \caption{Copied from \cite{Ch22} Figure 3. The first two are of type ${\rm II}_1$, the third is of type ${\rm III}_1$, and the last is of type ${\rm I}_1$.}\label{fig:piece}
\end{figure}

All the other symbols can also be chopped, as ensured by the identities given in the figures in Section \ref{sec:pictures}.
\end{proof}

Let $\mathsf{S}_u(F)$ denote the set of formal sums $\sum_s\mathfrak{a}_sC_s$ satisfying Lemma \ref{lem:substitution} (i);
let $\mathsf{S}_d(F)$ denote the set of formal sums $\sum_tD_t\mathfrak{b}_t$ satisfying Lemma \ref{lem:substitution} (ii).
Unlike in \cite{Ch22}, now it is no longer likely to find a preferred element of $\mathsf{S}_u(F)$ or $\mathsf{S}_d(F)$, so we need to deal with all possible choices.

For a generic $\mathsf{s}=\sum_s\mathfrak{a}_sC_s\in\mathsf{S}_u(F^\sharp)$ in the sense that each $(K^F|F^\sharp|C_s)$ is in generic position,
put
\begin{align}
\mathfrak{s}^\sharp(K,F,\mathsf{s})=U(K,F)+q^{\hat{\epsilon}(K,F)}{\sum}_s\mathfrak{a}_s(K^F|F^\sharp|C_s)\in\mathcal{F};
\end{align}
similarly, for a generic $\mathsf{t}=\sum_tD_t\mathfrak{b}_t\in\mathsf{S}_d(F^\flat)$ in the sense that each $(K_F|F^\flat|D_t)$ is in generic position, put
\begin{align}
\mathfrak{s}_\flat(K,F,\mathsf{t})=D(K,F)+q^{\check{\epsilon}(K,F)}{\sum}_t(K_F|F_\flat|D_t)\mathfrak{b}_t\in\mathcal{F}.
\end{align}

For a $1$-submanifold $X\subset\Sigma\times[0,1]$, let $\mathfrak{A}_3(X)$ denote the set of degree $3$ arcs of $X$.

For each knot $K$, recursively define the set $\mathcal{A}(K)\subset\mathcal{T}$ of {\it admissible expressions}.
When $|K|\le 3$, put $\mathcal{A}(K)=\{\Theta(K)\}$.
Suppose $|K|>3$ and that $\mathcal{A}(J)$ has been defined for each knot $J$ simpler than $K$. Take $F\in\mathfrak{A}_3(K)$,
take generic $\mathsf{s}\in\mathsf{S}_u(F)$, $\mathsf{t}\in\mathsf{S}_d(F)$,
and write $\mathfrak{s}^\sharp(K,F,\mathsf{s})=\sum_i\mathfrak{c}_iM_i$, $\mathfrak{s}_\flat(K,F,\mathsf{t})=\sum_jN_j\mathfrak{d}_j$, with $\mathfrak{c}_i,\mathfrak{d}_j\in\mathcal{T}$ and $M_i,N_j$ being knots subsequent to $K$;
put
\begin{align}
\mathcal{A}_u(K,F,\mathsf{s})&=\Big\{{\sum}_i\mathfrak{c}_i\mathfrak{g}_i\colon \mathfrak{g}_i\in\mathcal{A}(M_i)\Big\},  \label{eq:Au}  \\
\mathcal{A}_d(K,F,\mathsf{t})&=\Big\{{\sum}_j\mathfrak{h}_j\mathfrak{d}_j\colon \mathfrak{h}_j\in\mathcal{A}(N_j)\Big\}.  \label{eq:Ad}
\end{align}
Set
$$\mathcal{A}(K)={\bigcup}_{F\in\mathfrak{A}_3(K)}
\Big(\Big({\bigcup}_{\mathsf{s}}\mathcal{A}_u(K,F,\mathsf{s})\Big)\bigcup\Big({\bigcup}_{\mathsf{t}}\mathcal{A}_d(K,F,\mathsf{t})\Big)\Big).$$

For each arc $A\in\mathcal{H}_0(\mathsf{x},\mathsf{y})$ with $\mathsf{x},\mathsf{y}\in\Sigma\times\{0\}$, recursively define the set
$\mathcal{A}_{\rm pu}(A)\subset\mathcal{T}\mathcal{S}(\mathsf{x},\mathsf{y})$ of {\it purely upward admissible expressions}.
Put $\mathcal{A}_{\rm pu}(A)=\{\Theta(A)\}$ if $|A|\le2$; put $\mathcal{A}_{\rm pu}(A)=\{\mathfrak{s}_u(A)\}$ if $|A|=3$. Suppose $|A|>3$ and that $\mathcal{A}_{\rm pu}(B)$ has been defined for each arc $B$ simpler than $A$. Take $F\in\mathfrak{A}_3(A)$, and define $A^F$, $U(A,F)$, ``arcs subsequent to $A$", and so forth, similarly as above.
Take a generic $\mathsf{s}=\sum_s\mathfrak{a}_sC_s\in\mathsf{S}_u(F^\sharp)$.
Writing
\begin{align*}
s^\sharp(A,F,\mathsf{s}):=U(A,F)+q^{\hat{\epsilon}(A,F)}{\sum}_s\mathfrak{a}_s(A^F|F^\sharp|C_s)
\end{align*}
as ${\sum}_i\mathfrak{c}_iG_i$, where the $G_i$'s are arcs subsequent to $A$, put
\begin{align*}
\mathcal{A}_u(A,F,\mathsf{s})=\Big\{{\sum}_i\mathfrak{c}_i\mathfrak{g}_i\colon \mathfrak{g}_i\in\mathcal{A}_{\rm pu}(G_i)\Big\}.
\end{align*}
Set
$$\mathcal{A}_{\rm pu}(A)={\bigcup}_{F\in\mathfrak{A}_3(A)}\Big({\bigcup}_{\mathsf{s}}\mathcal{A}_u(A,F,\mathsf{s})\Big).$$

For each arc $A\in\mathcal{H}_1(\mathsf{x},\mathsf{y})$ with $\mathsf{x},\mathsf{y}\in\Sigma\times\{1\}$, define the set
$\mathcal{A}_{\rm pd}(A)\subset\mathcal{S}(\mathsf{x},\mathsf{y})\mathcal{T}$ of {\it purely downward admissible expressions} in a parallel way.

\begin{lem}
Suppose $K={\rm tr}(AB)$, where $A,B$ are arcs with ${\rm Cr}(B/A)=\emptyset$.
If $\sum_i\mathfrak{a}_i[C_i]\in\mathcal{A}_{\rm pu}(A)$, $\sum_j[D_j]\mathfrak{b}_j\in\mathcal{A}_{\rm pd}(B)$, and $\mathfrak{f}_{i,j}\in\mathcal{A}({\rm tr}(C_iD_j))$ for all $i,j$, then
${\sum}_{i,j}\mathfrak{a}_i\mathfrak{f}_{i,j}\mathfrak{b}_j\in\mathcal{A}(K).$
\end{lem}

For a stacked link $L=K_1\cdots K_r$, put
$$\mathcal{A}(L)=\big\{\mathfrak{e}_1\cdots\mathfrak{e}_r\colon\mathfrak{e}_1\in\mathcal{A}(K_1),\ldots,\mathfrak{e}_r\in\mathcal{A}(K_r)\big\}.$$
Suppose $M=S_1\sqcup\cdots\sqcup S_r$ is a multicurve. Each $\sigma\in{\rm Sym}(r)$ (the permutation group of $\{1,\ldots,r\}$) gives rise to a stacked link $S_{\sigma(1)}\cdots S_{\sigma(r)}$, which is not equivalent to $S_1\cdots S_r$ unless $\sigma$ is the identity.
Put
\begin{align}
\mathcal{A}(M)={\bigcup}_{\sigma\in{\rm Sym}(r)}\mathcal{A}(S_{\sigma(1)}\cdots S_{\sigma(r)}).   \label{eq:multi-curve}
\end{align}

The construction based on Lemma \ref{lem:substitution} has ensured $\mathcal{A}(K)\ne\emptyset$ for any knot $K$.
Hence $\mathcal{A}(M)\ne\emptyset$ for any multicurve $M$, establishing the surjectivity of $\theta$.

\section{Relations can be localized}

For a $1$-submanifold $X$, let $\Sigma(X)=\Sigma(v_1,\ldots,v_r)$ if $X\cap\Gamma_v\ne\emptyset$ exactly for $v=v_1,\ldots,v_r$.


Let $\Lambda$ be the set of $(v_1,\ldots,v_k)$ with $1\le v_1<\cdots<v_k\le n$ and $2\le k\le 6$.
For each $\vec{v}=(v_1,\ldots,v_k)\in\Lambda$, let $\Sigma(\vec{v})=\Sigma(v_1,\ldots,v_k)$,
and let $\mathcal{Z}(\vec{v})\subset\mathcal{T}$ denote the $R$-module of polynomials $\mathfrak{u}$ in $t_{i_1\cdots i_r}$'s with $r\le 3$ and $i_1,\ldots,i_r\in\{v_1,\ldots,v_k\}$ such that $|\mathfrak{u}|\le 6$ and $\mathfrak{u}=0$ in $\mathcal{S}(\Sigma(\vec{v}))$.
For $\mathfrak{u}\in\mathcal{Z}(\vec{v})$, we say that the relation $\mathfrak{u}=0$ is {\it supported} by $\Sigma(v_1,\ldots,v_k)$.

Let $\mathcal{I}$ denote the two-sided ideal of $\mathcal{T}$ generated by $\bigcup_{\vec{v}\in\Lambda}\mathcal{Z}(\vec{v})$.

\begin{thm}\label{thm:relation}
The ideal of defining relations of $\mathcal{S}$ is $\mathcal{I}$, i.e., $\mathcal{I}=\ker\theta$.
\end{thm}

Let $\zeta:\mathcal{T}\twoheadrightarrow\mathcal{T}/\mathcal{I}$ denote the quotient map.
For $\mathfrak{e},\mathfrak{e}'\in\mathcal{T}$, we say that $\mathfrak{e}$ is {\it congruent} to $\mathfrak{e}'$ and denote $\mathfrak{e}\equiv\mathfrak{e}'$, if $\mathfrak{e}-\mathfrak{e}'\in\mathcal{I}$, i.e. $\zeta(\mathfrak{e})=\zeta(\mathfrak{e}')$.

For a knot $K$, if elements of $\mathcal{A}(K)$ are congruent to each other, 
then we denote $\check{\mathfrak{a}}(K)\in\mathcal{T}/\mathcal{I}$ for the unique element of $\zeta(\mathcal{A}(K))$ and say that
{\it $\check{\mathfrak{a}}(K)$ is well-defined}.
Note that $\check{\mathfrak{a}}(K)$ is well-defined whenever $|K|\le 6$.

Given $\Omega=\sum_ia_iN_i$, with $a_i\in R$ and $N_i$ being a stacked link (probably a knot), let $\mathcal{A}(\Omega)=\{\sum_ia_i\mathfrak{e}_i\colon \mathfrak{e}_i\in\mathcal{A}(N_i)\}$; say that $\check{\mathfrak{a}}(\Omega)$ is well-defined if $\#\zeta(\mathcal{A}(\Omega))=1$. For two linear combinations $\Omega_1,\Omega_2$, denote $\Omega_1\equiv \Omega_2$ if $\check{\mathfrak{a}}(\Omega_1)$, $\check{\mathfrak{a}}(\Omega_2)$ are well-defined and equal, i.e., $\mathfrak{e}_1\equiv\mathfrak{e}_2$ for any $\mathfrak{e}_1\in\mathcal{A}(\Omega_1)$, $\mathfrak{e}_2\in\mathcal{A}(\Omega_2)$.

\begin{figure}[H]
  \centering
  \includegraphics[width=11.5cm]{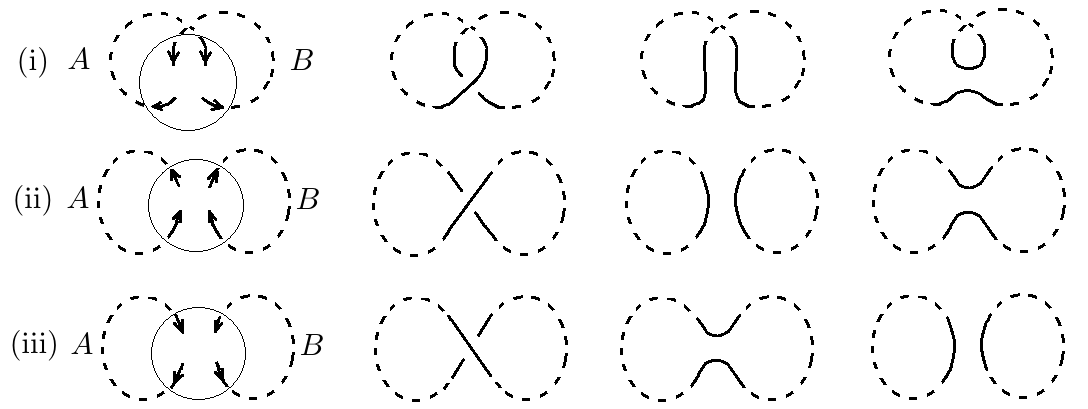}\\
  \caption{First row: the EP $(A,B)$; $L_\times={\rm tr}(A){\rm tr}(B)$, $L_\infty={\rm tr}(AB)$, $L_0={\rm tr}(A\overline{B})$.
  Second row: the EP $(A,B)$; $L_\times={\rm tr}(AB)$, $L_\infty={\rm tr}(A){\rm tr}(B)$, $L_0={\rm tr}(A\overline{B})$.
  Third row: the EP $(A,B)$; $L_\times={\rm tr}(AB)$, $L_\infty={\rm tr}(A\overline{B})$, $L_0={\rm tr}(A){\rm tr}(B)$.}\label{fig:elementary}
\end{figure}

An {\it elementary skein triple} (EST) is a triple of the form $(L_\times,L_\infty,L_0)$ in one of the three cases shown in Figure \ref{fig:elementary}.

Here is the main procedure of establishing Theorem \ref{thm:relation}. It will be shown that $\check{\mathfrak{a}}(K)$ is well-defined for each knot $K$, so that $\check{\mathfrak{a}}(L)$ is well-defined for each stacked link $L$. The transformation from $L$ to $\Theta(L)$ can be implemented via crossing-resolving of specific kind in the following sense: there exists a series of linear-combinations of stacked links $L=\Omega_0,\Omega_1,\ldots,\Omega_m=\Theta(L)$ such that
$$\Omega_{i+1}-\Omega_i=\mathfrak{a}\big(L_\times-q^{\frac{1}{2}}L_\infty-\overline{q}^{\frac{1}{2}}L_0\big)\mathfrak{b}$$
for some EST $(L_\times,L_\infty,L_0)$ and $\mathfrak{a},\mathfrak{b}\in\mathcal{T}$.
We will show that $L_\times\equiv q^{\frac{1}{2}}L_\infty+\overline{q}^{\frac{1}{2}}L_0$ always. Thus, $L\equiv\Theta(L)$. In particular, $\mathfrak{m}\equiv\Theta(\mathfrak{m})$ for each monomial $\mathfrak{m}$. As a consequence, whenever $\mathfrak{f}=0$ in $\mathcal{S}$, we have $\mathfrak{f}\equiv\Theta(\mathfrak{f})=0$, i.e. $\mathfrak{f}\in\mathcal{I}$.

Technically, we need to prove several intermediate lemmas by ``tortuous" inductions.

For $m\ge6$ and $c\ge 0$, let $\phi_{m,c}$ stand for the statement that $\check{\mathfrak{a}}(K)$ is well-defined for each knot $K$ with $\lambda(K)=(m,c)$; let $\Phi_{m,c}$ stand for the statement that $\phi_{m',c'}$ holds for all $(m',c')\preceq(m,c)$.

\begin{lem}
Suppose $\Phi_{m,c}$ holds, and $K$ is a knot with $\lambda(K)=(m,c)$.
\begin{enumerate}
  \item[\rm(a)] If $F\in\mathfrak{A}_3(K)$, then $K\equiv\mathfrak{s}^\sharp(K,F,\mathsf{s})\equiv\mathfrak{s}_\flat(K,F,\mathsf{t})$
       for all generic $\mathsf{s}\in\mathsf{S}_u(F)$ and $\mathsf{t}\in\mathsf{S}_d(F)$.
  \item[\rm(b)] For any arc $F\subset K$,
       $$K\equiv q^{\hat{\epsilon}(K,F)}K^F+U(K,F)\equiv q^{\check{\epsilon}(K,F)}K_F+D(K,F).$$
\end{enumerate}
\end{lem}

\begin{lem}
Suppose $\Phi_{m,c}$ holds. Then
\begin{enumerate}
  \item[\rm(a)] $L_\times\equiv q^{\frac{1}{2}}L_\infty+\overline{q}^{\frac{1}{2}}L_0$
        for each EST $(L_\times,L_\infty,L_0)$ with $\lambda(L_\times)\preceq(m,c)$;
  \item[\rm(b)] $K_1K_2\equiv K_2K_1$ for any disjoint knots $K_1,K_2$ such that $\lambda(K_1K_2)\prec(m,c)$.
\end{enumerate}
\end{lem}

\begin{lem}
Suppose $\Phi_{m,0}$ holds.
\begin{enumerate}
  \item[\rm(a)] Let $S,S'$ be simple curves such that $|S|\le m$ and $S'$ results from shrinking a degree $2$ arc of $S$. Then $S\equiv S'$.
  \item[\rm(b)] Let $S,S'$ be simple curves such that $|S|,|S'|\le m$ and $S'$ results from sliding $S$ across an interior vertex.
  Then $S\equiv S'$.
\end{enumerate}
\end{lem}

\begin{proof}[Proof for {\rm(b)}]
Assume $S={\rm tr}(AF)$, $S'={\rm tr}(AF')$, for some arcs $A,F,F'$ such that $|F|,|F'|\le4$ and ${\rm tr}(F'\overline{F})$ is a simple curve encircling an interior vertex.

Take $\sum_i\mathfrak{a}_i[C_i]\in\mathcal{A}_{\rm pu}(A)$.
For each $i$, we have ${\rm tr}(C_iF')\equiv{\rm tr}(C_iF)$, since ${\rm tr}(C_iF')={\rm tr}(C_iF)$ is a relation of degree at most $6$ supported by $\Sigma(C_i\cup F\cup F')$. Hence
$$\check{\mathfrak{a}}(S)-\check{\mathfrak{a}}(S')={\sum}_i\mathfrak{a}_i\big(\check{\mathfrak{a}}(C_iF)-\check{\mathfrak{a}}(C_iF')\big)=0.$$
\end{proof}

\begin{lem}
Suppose $\Phi_{m,1}$ holds.
\begin{enumerate}
  \item[\rm(a)] If $S_1,\ldots,S_r$ are disjoint simple curves with $|S_1|+\cdots+|S_r|\le m$, then
        $S_{\sigma(1)}\cdots S_{\sigma(r)}\equiv S_1\cdots S_r$ for each $\sigma\in{\rm Sym}(r)$.
        Consequently, elements of $\mathcal{A}(S_1\sqcup\cdots\sqcup S_r)$ are congruent to each other.
  \item[\rm(b)] If $M,M'$ are isotopic multicurves of degree at most $m$, then $M\equiv M'$.
\end{enumerate}
\end{lem}

\begin{proof}[Proof for {\rm(b)}]
This is based on Lemma \ref{lem:isotopy}.
\end{proof}

\begin{lem}
Suppose $\Phi_{m,c}$ holds. Then the following statements are true.
\begin{enumerate}
  \item[\rm(a)] $K\equiv \Theta(K)$ for any knot $K$ with $\lambda(K)=(m,c)$.
  \item[\rm(b)] $K_1K_2\equiv \Theta(K_1K_2)$ for any knots $K_1,K_2$ with $\lambda(K_1K_2)\preceq(m,c)$.
  \item[\rm(c)] If $\sum_ja_jK_j=0$ in $\mathcal{S}$, where each $K_j$ is a knot with $\lambda(K_j)\preceq(m,c)$, then $\sum_ja_jK_j\equiv 0$.
\end{enumerate}
\end{lem}

\begin{lem}
$\phi_{m,c}$ holds for all $m\ge 6$ and all $c\ge 0$.
\end{lem}

\begin{lem}\label{lem:link}
$K_1\cdots K_r\equiv\Theta(K_1\cdots K_r)$, for any knots $K_1,\ldots,K_r$.
\end{lem}

\medskip

\begin{proof}[Proof of Theorem \ref{thm:relation}]
Suppose $\mathfrak{f}=\sum_ia_i\mathfrak{g}_i=0$ in $\mathcal{S}$, where $a_i\in R$ and $\mathfrak{g}_i$ is a monomial. Since each $\mathfrak{g}_i$ itself is a stacked link, by Lemma \ref{lem:link}, $\mathfrak{g}_i\equiv \Theta(\mathfrak{g}_i)$.
Hence
$\mathfrak{f}\equiv{\sum}_ia_i\Theta(\mathfrak{g}_i)=\Theta(\mathfrak{f})=0.$
\end{proof}



\section{Examples of cutting systems}

\subsection{$\Sigma_{0,4}$ revised}

\begin{figure}[H]
  \centering
  \includegraphics[width=12cm]{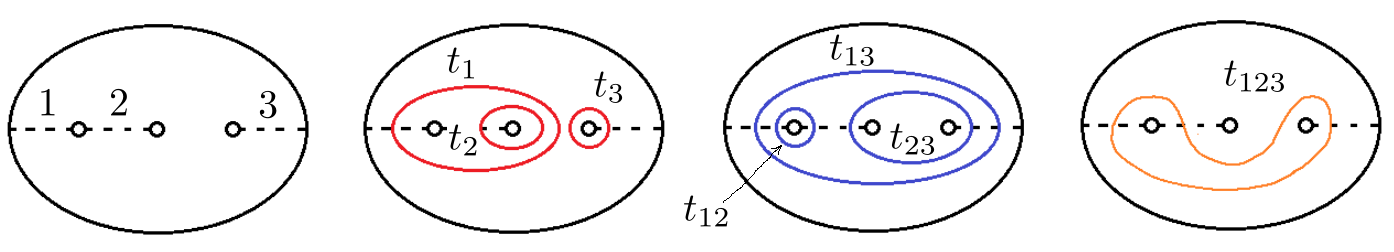}\\
  \caption{From left to right: a cutting system; $t_1,t_2,t_3$; $t_{12},t_{13},t_{23}$; $t_{123}$.}\label{fig:g0-new}
\end{figure}

Figure \ref{fig:g0-new} shows a cutting system, together with an associated generating set.

\subsection{$\Sigma_{1,3}$}

\begin{figure}[h]
  \centering
  \includegraphics[width=11cm]{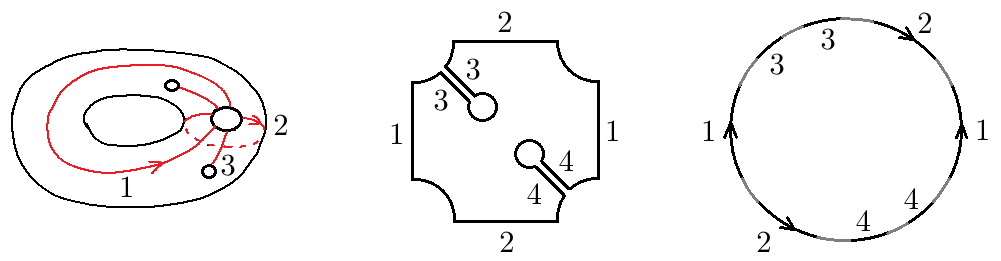}\\
  \caption{Left: $\Sigma_{1,3}$, with a cutting system $G$ drawn in red. Middle and right: the closure of $\Sigma_{1,3}\setminus G$.}\label{fig:Sigma1-3}
\end{figure}

\begin{figure}[H]
  \centering
  \includegraphics[width=12.5cm]{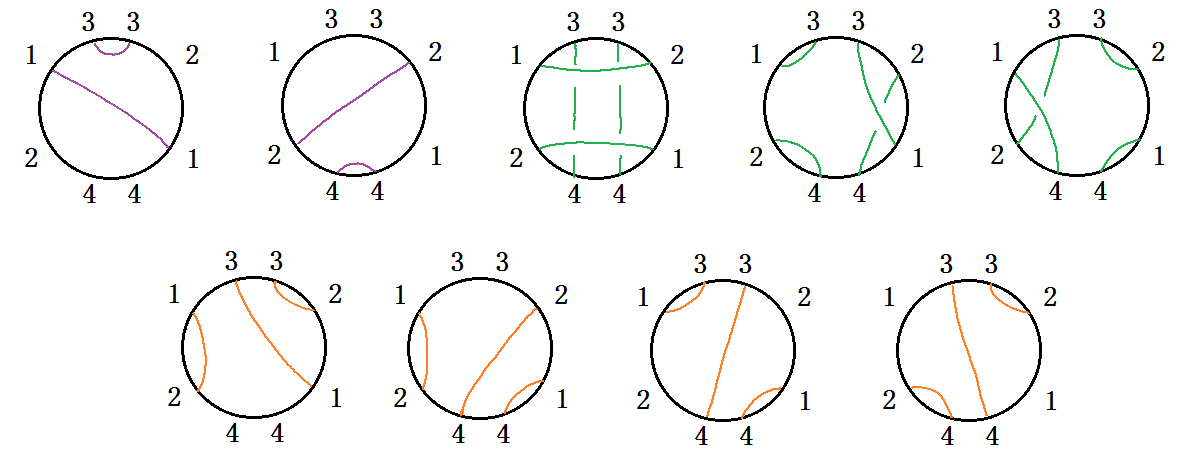}\\
  \caption{The $t_i$'s, $t_{ij}$'s, $t_{ijk}$'s are respectively drawn in purple, green, orange.}\label{fig:genus2-generator}
\end{figure}

\subsection{$\Sigma_{g,k+1}$}

\begin{figure}[h]
  \centering
  \includegraphics[width=12cm]{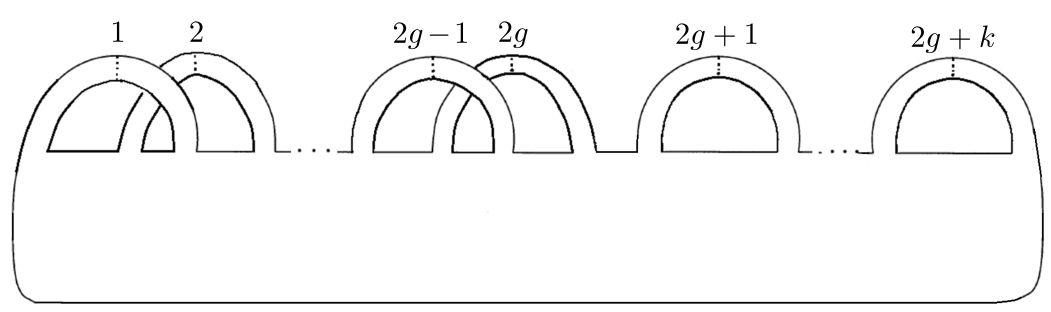}\\
  \caption{$\Sigma_{g,k+1}$; the dotted lines (from left to right) are $\gamma_1,\ldots,\gamma_{2g+k}$.}\label{fig:Sigma}
\end{figure}

Display $\Sigma_{g,k+1}$ as in Figure \ref{fig:Sigma}. A preferred cutting system is $\bigcup_{i=1}^{2g+k}\gamma_i$, where $\gamma_1,\ldots,\gamma_{2g+k}$ are indicated by the dotted lines.

\begin{figure}[H]
  \centering
  \includegraphics[width=11.5cm]{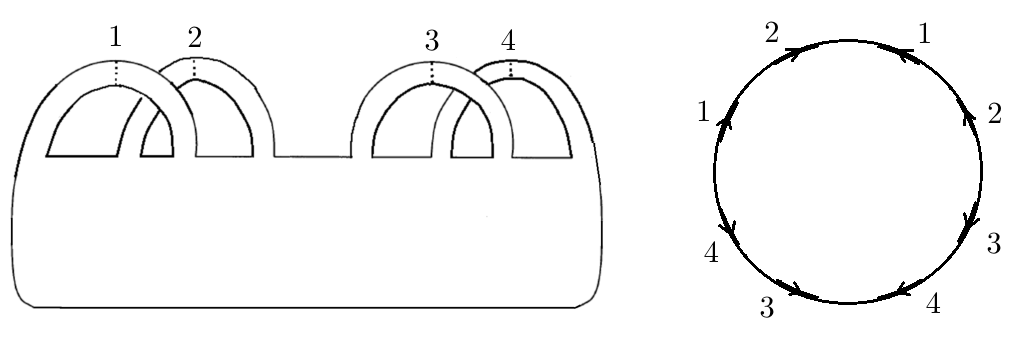}\\
  \caption{Left: $\Sigma_{2,1}$. Right: the closure of $\Sigma_{2,1}\setminus(\gamma_1\cup\cdots\cup\gamma_4)$.}\label{fig:Sigma2-1}
\end{figure}

In particular, see Figure \ref{fig:Sigma2-1} for $\Sigma_{2,1}$.

\begin{figure}[H]
  \centering
  \includegraphics[width=12.5cm]{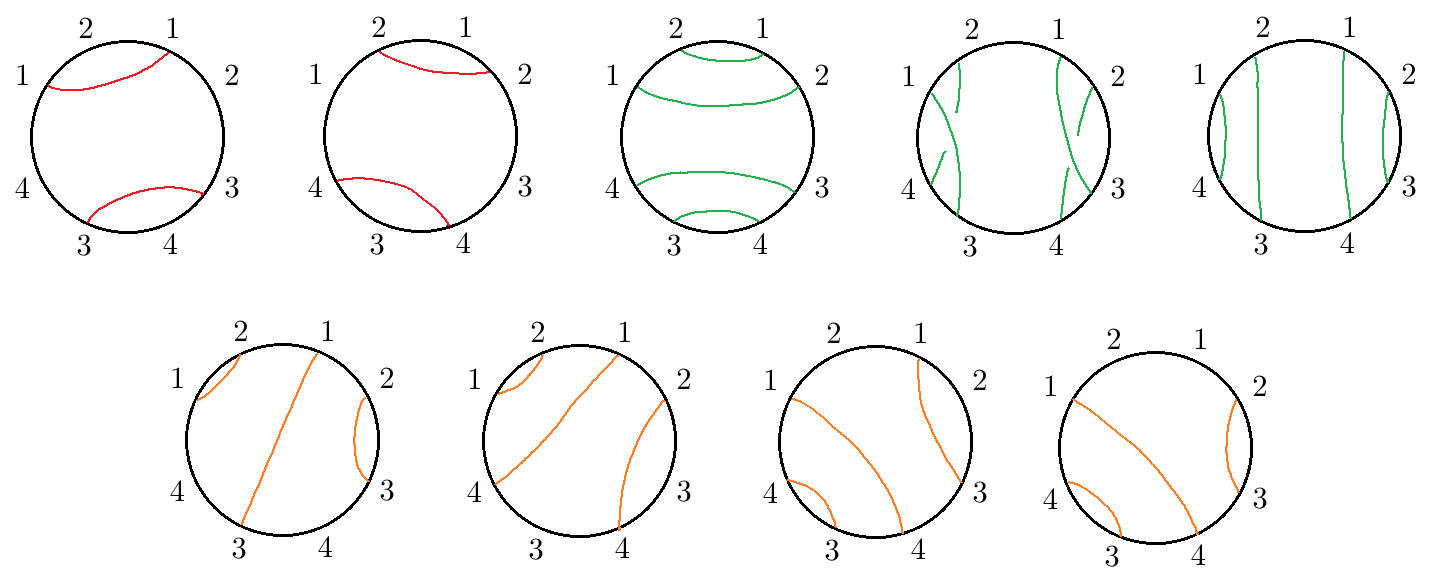}\\
  \caption{The $t_i$'s, $t_{ij}$'s, $t_{ijk}$'s are respectively drawn in red, green, orange.}\label{fig:genus2-generator}
\end{figure}

A generating set for $\mathcal{S}(\Sigma_{2,1})$ is shown in Figure \ref{fig:genus2-generator}.

A presentation for $\mathcal{S}(\Sigma_2)$ will be obtained by adding relations given by capping off the boundary.



\newpage

\section{Appendix: pictorial proofs for basic identities} \label{sec:pictures}

An apology: at first sight these formulas seem to be cluttered arbitrarily; but actually they are obtained based on various attempts, and have been arranged to be nearly optimal.

\begin{figure}[H]
  \centering
  \includegraphics[width=10.5cm]{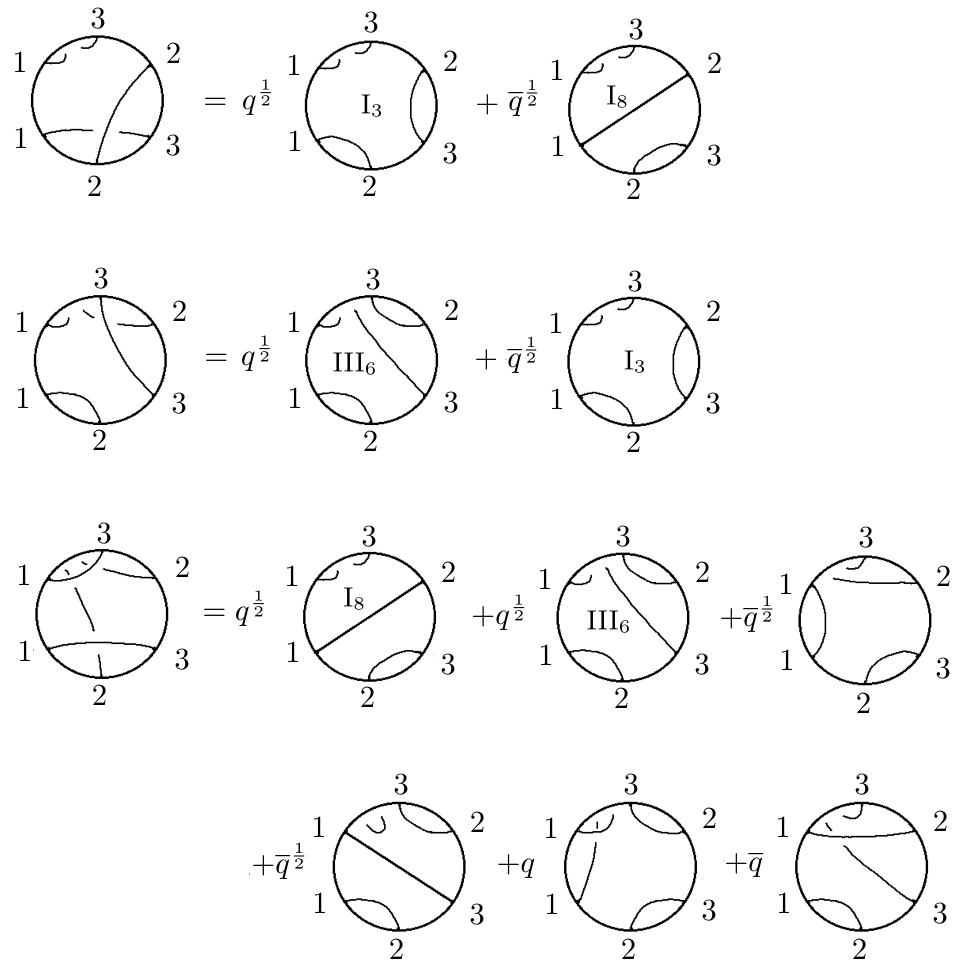}\\
  \caption{Taking charge of ${\rm I}_3$, ${\rm I}_8$, ${\rm III}_6$.}\label{fig:chop-3}
\end{figure}

\begin{figure}[H]
  \centering
  \includegraphics[width=12.5cm]{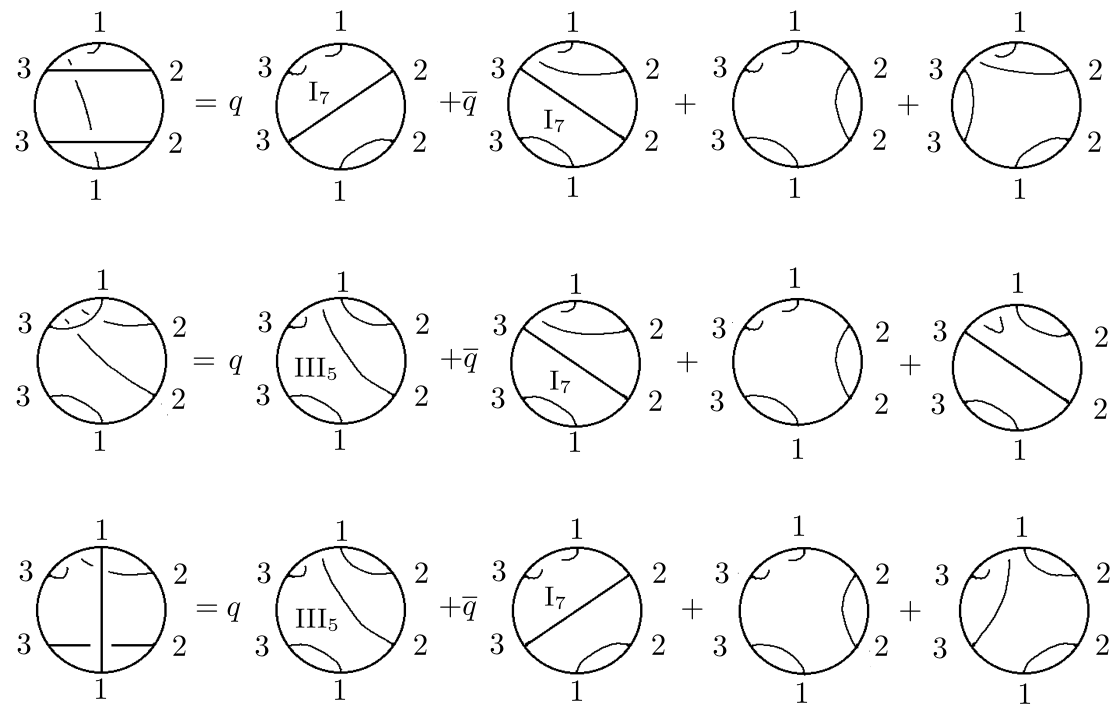}\\
  \caption{Taking charge of ${\rm I}_7$ and ${\rm III}_5$.}\label{fig:chop-1}
\end{figure}

\begin{figure}[H]
  \centering
  \includegraphics[width=12.5cm]{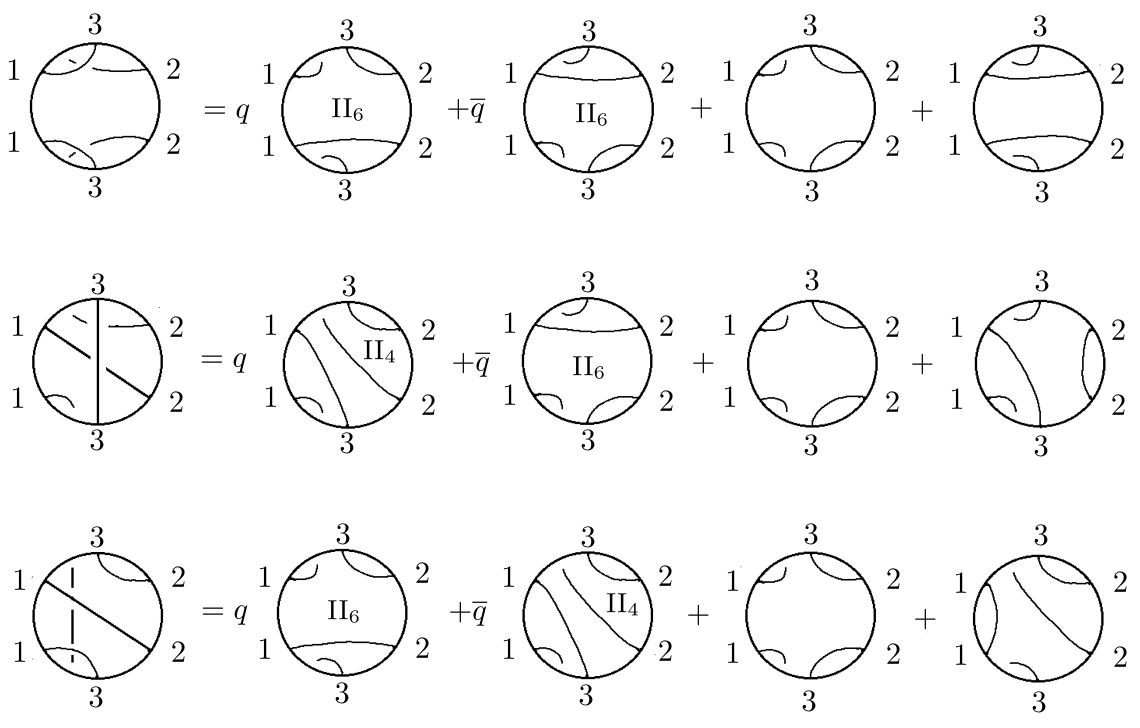}\\
  \caption{Taking charge of ${\rm II}_6$ and ${\rm II}_4$.}\label{fig:chop-2}
\end{figure}

\begin{figure}[H]
  \centering
  \includegraphics[width=8.5cm]{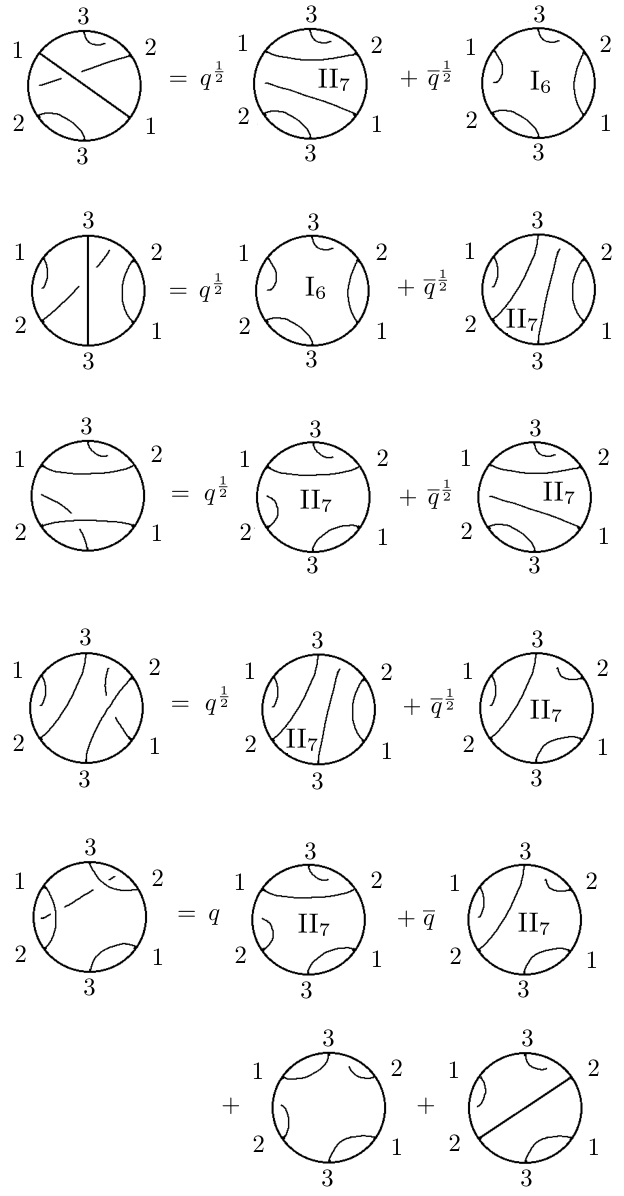}\\
  \caption{Taking charge of ${\rm I}_6$ and ${\rm II}_7$.}\label{fig:chop-4}
\end{figure}

\begin{figure}[H]
  \centering
  \includegraphics[width=8.5cm]{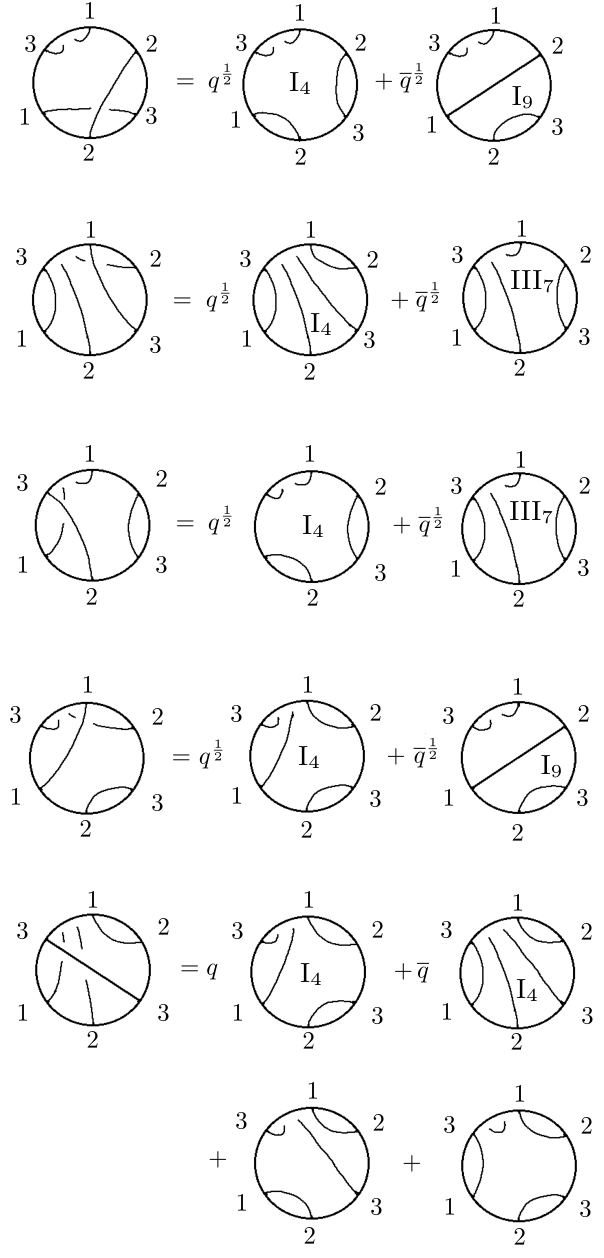}\\
  \caption{Taking charge of ${\rm I}_4$, ${\rm I}_9$ and ${\rm III}_7$.}\label{fig:chop-5}
\end{figure}

\begin{figure}[H]
  \centering
  \includegraphics[width=8cm]{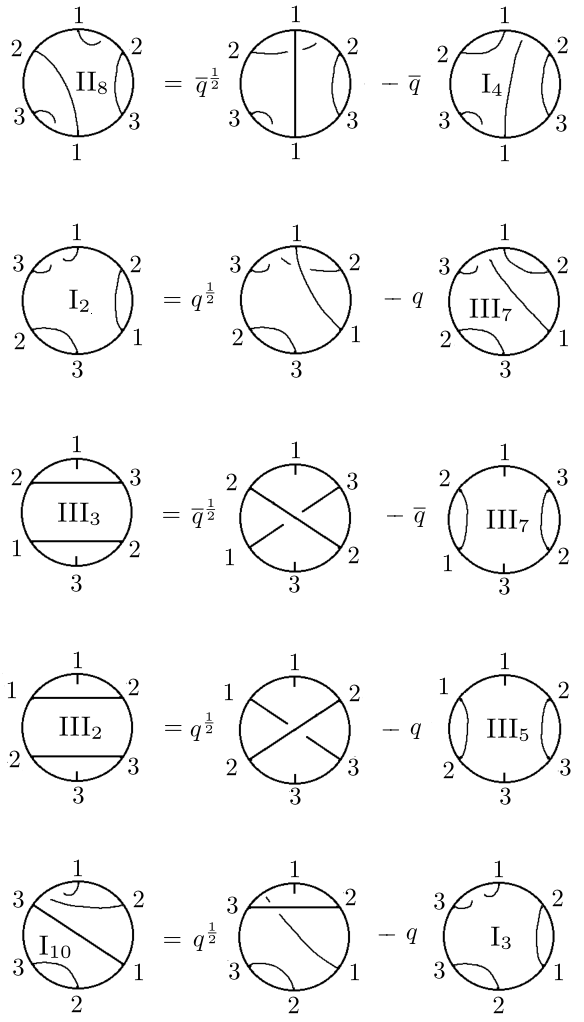}\\
  \caption{Taking charge of ${\rm II}_8$, ${\rm I}_2$, ${\rm III}_3$, ${\rm III}_2$, ${\rm I}_{10}$.}\label{fig:chop-6}
\end{figure}

\begin{figure}[H]
  \centering
  \includegraphics[width=8cm]{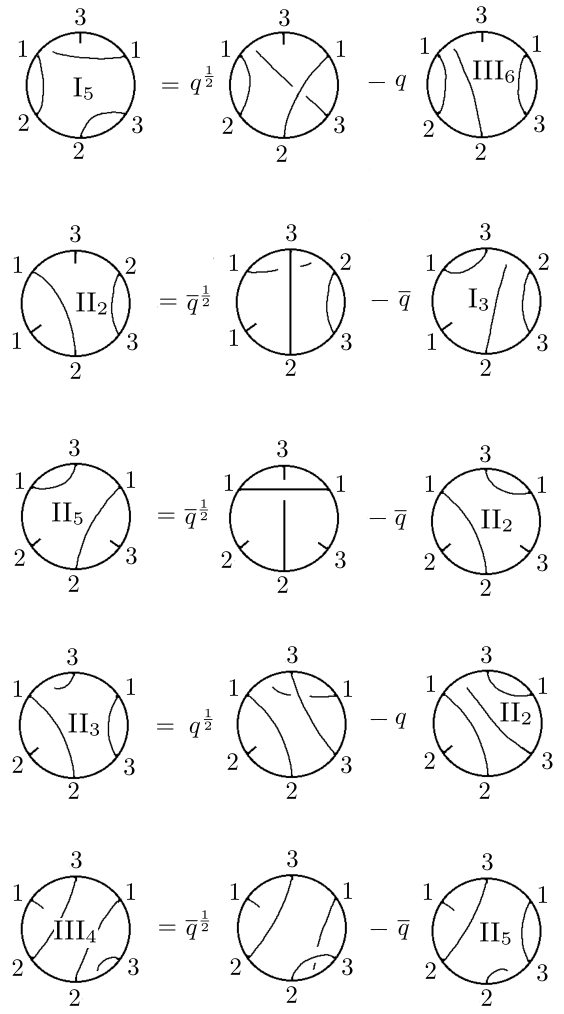}\\
  \caption{Taking charge of ${\rm I}_5$, ${\rm II}_2$, ${\rm II}_5$, ${\rm II}_3$, ${\rm III}_4$.}\label{fig:chop-7}
\end{figure}








\bigskip

\noindent
Haimiao Chen (orcid: 0000-0001-8194-1264)\ \ \  \emph{chenhm@math.pku.edu.cn} \\
Department of Mathematics, Beijing Technology and Business University, \\
Liangxiang Higher Education Park, Fangshan District, Beijing, China.

\end{document}